\documentclass{article}

\begin{document}
\begin{center}
{\large \bf  On the existence of coincidence and common fixed point of two rational type contractions and an application in dynamical programming  }\\
\end{center}
\vskip 0.3cm
\begin{center}
 {\small \sc Ahmed H. Soliman$^{a,c}$\footnote{mahmod@kku.edu.sa}  and  { TAMER NABIL$^{b,c}$\footnote{$t_{-}$3bdelsadek@yahoo.com} }
  \\
 \small{ $^{a}$Department of Mathematics, Faculty of Science, Al-Azhar University, Assiut 71524, Egypt}
 \\
 \small{ $^{b}$Suez Canal University, Faculty of Computers and Informatics, Department of Basic Science, Ismailia, Egypt}
 \\
\small{ $^{c}$King Khalid University, Faculty of Science, Mathematics Department, Abha 9004, Saudi Arabia}
}
\end{center}
\begin{abstract}
In this work, we establish some coincidence point results for self-mappings satisfying rational type contractions in generalized metric spaces  in the sense of Branciari \cite{Branciari}. Presented coincidence point theorems weak and extend numerous existing theorems in the literature besides
furnishing some illustrative examples for our results. Finally, our results applies, in particular, to the study of
 solvability of functional equations arising in dynamic programming.
\end{abstract}
\section{Introduction}
\quad Banach contraction principle is one of the most important aspects of fixed point theory as a source of the existence and uniqueness of solutions of many problems in various  branches inside and outside mathematics (see, \cite{Cesari,Hale,Rouche}). Some generalizations of this theorem replace the contraction condition by a weaker. For instance, in 1973, Dass and Gupta \cite{DassandGupta} defined the following rational type contraction which is more general than the contraction condition.
\begin{eqnarray}
d(Ax,Ay)\leq a d(x,y)+ \frac{b d(y,Ay)(d(x,Ax)+1)}{1+d(x,y)}\; \forall \; x,y \in X\; \textrm{and} \; a, b \geq 0,  \; a + b < 1.
\end{eqnarray}
where $A: X \to X $ be a mapping from a metric space $ X $ into itself .

Recently, 2015, Almeida, Rold´an-L´opez-de-Hierro and Sadarangani\cite{ARS2015} introduced an extension of the condition (1) of Dass and Gupta\cite{DassandGupta} as follows:
\begin{eqnarray}
d(Ax,Ay)\leq \phi(M(x,y))+ C\min\{d(x,Ax), d(y,Ay),d(x,Ay),d(y,Ax)\}\; \forall \; x,y \in X, \; C \geq 0,
\end{eqnarray}
where $M(x,y)$ is defined by
$$
M(x,y)=\max\{d(x,y),\frac{ d(x,Ax)(d(y,Ay)+1)}{1+d(x,y)}, \frac{ d(y,Ay)(d(x,Ax)+1)}{1+d(x,y)}\},
$$
and $\phi:[0,\infty) \to [0,\infty)$ is a non-decreasing upper semi-continuous function with
$\phi(t)< t $ for all $t > 0.$

Another offshoot of generalizations of Banach's theorem is based on extending the axioms of metric spaces. It worth to mention that the use of triangle inequality in a metric space $(X,d)$ is of extreme importance since it implies that $d$ is continuous, each open ball is an open set, a sequence may converge to unique point, every convergent sequence is a Cauchy sequence and other things. In 2000, Branciari \cite{Branciari} introduced a new concept of generalized metric space by replacing the triangle inequality of a metric space by a so-called rectangular inequality. Since then, various works have dealt with fixed point results in such spaces (see, \cite{Arshad2013,Asadi,Azam2008,Das2007,Das2009,Fora}). It was not directly noted that such generalized metric spaces (G.M.S, for short) may fail to satisfy the conditions which mentioned above in metric spaces.

In this paper, we introduce coincidence point theorems for two contraction self-mappings of rational type in generalized metric spaces. Our result improve the results due to  Almeida, Rold´an-L´opez-de-Hierro and Sadarangani \cite{ARS2015}. These theoretical theorems are applied to the study of the existence solutions to a system of functional equations in dynamic programming.
\section{Preliminaries}
\quad In this section, we present some preliminaries and notations related to  rational type contraction and generalized metric spaces.
\\
{\bf Definition 2.1} (Branciari \cite{Branciari}). Suppose that $X$ be a nonempty set and $d: X\times X \to [0,\infty)$ be a distance function such that for all $w, x, y,z  \in X$ and $w\ne x\ne y\ne z,$

(i)  $d(w,x)=0$ $\Leftrightarrow$ $w = x,$

(ii) $d(w,x)=d(x,w),$

(iii) $d(w,x)\leq d(x,y)+d(y,z)+d(z,w)$ (quadrilateral inequality).
\\
Then we called $(X,d)$ G.M.S.

The following example show that G.M.S more general than metric spaces
 \\
 {\bf Example 2.1. } Suppose that $X= \{\frac{5}{6}, \frac{2}{3},\frac{7}{12},\frac{8}{15}\}.$  Define $d$ on $ X\times X $ as follows:
 \begin{eqnarray*}
 && d(\frac{5}{6},\frac{2}{3})=d(\frac{7}{12},\frac{8}{15})=\frac{4}{9}, \;\; d(\frac{5}{6},\frac{8}{12})=d(\frac{2}{3},\frac{7}{12})=\frac{1}{3},
 \\&& d(\frac{5}{6},\frac{7}{12})=d(\frac{2}{3},\frac{8}{12})=\frac{8}{9}, \;\; d(x,x)=0,\;  d(x,y)=d(y,x).
 \end{eqnarray*}
 Then $(X, d)$ is a G.M.S but not metric space.
 \\
 {\bf Remark 2.1.} We note that Definition 2.1 especially condition (iii) dose not ensure that $d$ is continuous on its domain, see \cite{Branciari}. However, Convergent sequences,  Cauchy sequences and  completeness in G.M.S  may not hold.
  \\
 {\bf Definition 2.2} ( Rosa and Vetro\cite{RV2014}). Suppose that $(X,d)$ be a G.M.S and let $\{x_{n}\}$ be a sequence in $X.$ Then

 (i) $\{x_{n}\}$ converges to $x\in X $ in G.M.S iff $\lim\limits_{n\to \infty} d(x_{n},x)=0,$

 (ii) $\{x_{n}\}$ is a Cauchy in G.M.S iff $\forall\; \epsilon >0, \;  \exists \; K(\epsilon)>0 $ such that  $ d(x_{r},x_{s})<\epsilon,\; \forall \; r > s\geq K(\epsilon),$

 (iii) $(X,d)$ is called complete G.M.S if every Cauchy sequence in $X$ converges to a point in $X.$

 In 2009, Sarma et al. \cite{Sarma2009} introduced the following example which show Remark 2.1.
 \\
 {\bf Example 2.2} (Sarma et al. \cite{Sarma2009}). Suppose that $X = D\cup E,$ where $D=\{0,2\}$ and $E = \{\frac{1}{n}: n\in N (\textrm{ the set of all natural numbers } )\}.$ Define $ d $ from $ X\times X $ into $ [0,+\infty) $ as follows:
 $$
d(u,v)= \left \{
\begin{array}{lll}
 0, & u = v \\
 1, & u\ne v \; \&\; \{u,v\}\subset D \; \textrm{or}\; \{u,v\}\subset E,\\
\end{array} \right.
$$
and
$d(u,v) = d(v,u) = u$ if $u\in D $ and $v\in E.$
\\
Then $(X, d)$ is a complete G.M.S. Moreover, one can see that:

 (1) $d(\frac{1}{n}, 0)=0$ and $d(\frac{1}{n}, 2)=2$ $\Rightarrow$ $\{\frac{1}{n}\}$ is not Cauchy sequence.

 (2) $d(\frac{1}{n}, \frac{1}{2})\ne d(\frac{1}{2}, 0)$ $\Rightarrow$ $d$ is not continuous.
 \\
{\bf Definition 2.3} (\cite{RV2014}). Let $A, B : X \to X$ and $\beta : X \times X \to [0,\infty).$ The mapping $A$ is $B-\beta-$admissible if, for all $x,y\in X$ such that $\beta(Bx,By) > 1,$ we have $\beta(Ax,Ay) > 1.$ If $B$ is
the identity mapping, then $A$ is called $\beta-$admissible.
\\
{\bf Definition 2.4} (\cite{RV2014}). Let $(X, D)$ be a G.M.S and $\beta : X \times X \to [0,\infty).$  $X$ is $\beta-$regular if, for
each sequence $\{x_{n}\}$ in  $X$  such that $\beta(x_{n}, x_{n+1})> 1 $ for all $ n \in N $ and $ x_{n} \to x,$ then
there exists a subsequence $\{x_{n_{k}}\}$ of $\{ x_n \}$ such that $\beta(x_{n_k},x) > 1 $ $\forall$ $k \in N.$
\section{Main results}
\quad In this section we introduce some coincidence point results for two rational contraction self-mappings on
G.M.S.
\\
{\bf Theorem 3.1.} Let $(X,d)$ be a G.M.S and let $A$ and $B$ be self-mappings on $X$ such that $AX\subset BX.$  Suppose that $(BX,d)$ is a complete G.M.S and the following condition holds:
 \begin{eqnarray}
d(Ax,Ay)\leq \phi(M(x,y))+ C \min\{d(Bx,Ax), d(By,Ay),d(Bx,Ay),d(By,Ax)\}\; \forall \; x,y \in X, \; C \geq 0,
\end{eqnarray}
where $M(x,y)$ is defined by
$$
M(x,y)=\max\{d(Bx,By),\frac{ d(Bx,Ax)(d(By,Ay)+1)}{1+d(Bx,By)}, \frac{ d(By,Ay)(d(Bx,Ax)+1)}{1+d(Bx,By)}\}.
$$
and $\phi: [0,\infty) \to  [0,\infty) $ be a continuous, nondecreasing function and $\phi(t)=0 \; \iff t=0.$
\\
Then $A$ and $B$ have a unique point of coincidence in $X.$ Moreover if $A$ and $B$ are weakly compatible, then $A$ and $B$ have a unique common fixed point.
\\
{\bf Proof.} Define the sequence $\{x_n\}$ and $\{z_n\}$ in $X$ defined by
$$
z_n =Bx_{n+1}=Ax_{n}.
$$
If $z_{n}=z_{n+1},$ then $z_{n+1}$ is a point of coincidence   of $A$ and $B.$ Consequently, we can suppose that $z_{n}\ne z_{n+1} $ for all $n\in N.$
\\
Now, by (3), we have
\begin{eqnarray}
\nonumber & d(Ax_{n},Ax_{n+1})&\leq \phi(M(x_{n},x_{n+1}))+C \min\{d(Bx_{n},Ax_{n}),d(Bx_{n+1},Ax_{n+1}),d(Bx_{n},Ax_{n+1}),d(Bx_{n+1},Ax_{n})\}
\\&&= \phi(M(x_{n},x_{n+1}))
\end{eqnarray}
where
\begin{eqnarray}
\nonumber &M(x_{n},x_{n+1})&=\max\{d(Bx_{n},Bx_{n+1}),\frac{ d(Bx_{n},Ax_{n})(d(Bx_{n+1},Ax_{n+1})+1)}{1+d(Bx_{n},Bx_{n+1})}, \frac{ d(Bx_{n+1},Ax_{n+1})(d(Bx_{n},Ax_{n})+1)}{1+d(Bx_{n},Bx_{n+1})}\}
\\&&\nonumber = \max\{d(z_{n-1},z_{n}),\frac{d(z_{n-1},z_{n})(1+d(z_{n},z_{n+1}))}{1+d(z_{n-1},z_{n})}, d(z_{n},z_{n+1})\},
\end{eqnarray}
we consider the following cases

$\bullet$ If $M(x_{n},x_{n+1})= d(z_{n-1},z_{n})$ from $(4)$ we have
\begin{eqnarray}
 & d(z_{n},z_{n+1})&\leq \phi(d(z_{n-1},z_{n}))< d(z_{n-1},z_{n})
\end{eqnarray}

 $\bullet$ If $M(x_{n},x_{n+1})= \frac{d(z_{n-1},z_{n})(1+d(z_{n},z_{n+1}))}{1+d(z_{n-1},z_{n})}$ from $(4)$ we obtain
 \begin{eqnarray}
 \nonumber & d(z_{n},z_{n+1})&\leq \phi(\frac{d(z_{n-1},z_{n})(1+d(z_{n},z_{n+1}))}{1+d(z_{n-1},z_{n})})< \frac{d(z_{n-1},z_{n})(1+d(z_{n},z_{n+1}))}{1+d(z_{n-1},z_{n})}.
\end{eqnarray}
Hence
\begin{eqnarray}
\nonumber  & d(z_{n},z_{n+1})& < d(z_{n-1},z_{n}),
\end{eqnarray}
that is (5) holds.

$\bullet$ If $M(x_{n},x_{n+1})= d(z_{n},z_{n+1})$ from $(4)$ we get
\begin{eqnarray}
\nonumber  & d(z_{n},z_{n+1})& < d(z_{n},z_{n+1}),
\end{eqnarray}
which is impossible.

In any case, we proved that (5) holds. since $\{d(z_{n},z_{n+1})\}$ is decreasing sequence. Hence, it converges to a nonnegative number, $s\geq 0.$ If $s > 0,$ then letting $n\to +\infty$ in (4), we deduce
\begin{eqnarray}
\nonumber & s &\leq \phi(\max\{s,\frac{s(1+s)}{1+s}, s\})= \phi(s)< s,
\end{eqnarray}
which implies that $s = 0,$ that is
\begin{eqnarray}
& \lim\limits_{n\to\infty}d(z_{n},z_{n+1})& = 0.
\end{eqnarray}
Suppose that $z_{n}\ne z_{m}$ for all $m\ne n$ and prove that $\{z_{n}\}$ is G.M.S Cauchy sequence. First, we show that the sequence $\{d(z_{n},z_{n+2})\}$ is bounded. Since $\lim\limits_{n\to\infty}d(z_{n},z_{n+1})=0,$ there exists $L > 0 $ such that $d(z_{n},z_{n+1})\leq L $ for all $n\in N.$ If $d(z_{n},z_{n+2})> L $ for all $n\in N,$ from (3) we have
\begin{eqnarray}
\nonumber & d(z_{n},z_{n+2})&=d(Ax_{n},Ax_{n+2})
\\&&\nonumber \leq \phi(M(x_{n},x_{n+2}))+C \min\{d(Bx_{n},Ax_{n}),d(Bx_{n+2},Ax_{n+2}),d(Bx_{n},Ax_{n+2}),d(Bx_{n+2},Ax_{n})\}
\\&&\nonumber = \phi(M(x_{n},x_{n+2}))+ C \min\{d(z_{n-1},z_{n}),d(z_{n+1},z_{n+2}),d(z_{n-1},z_{n+2}),d(z_{n+1},z_{n})\}
\\&&= \phi(M(x_{n},x_{n+1}))\;\; as \;\; n\to \infty,
\end{eqnarray}
where
\begin{eqnarray}
\nonumber &M(x_{n},x_{n+2})&=\max\{d(Bx_{n},Bx_{n+2}),\frac{ d(Bx_{n},Ax_{n})(d(Bx_{n+2},Ax_{n+2})+1)}{1+d(Bx_{n},Bx_{n+2})},
\frac{ d(Bx_{n+2},Ax_{n+2})(d(Bx_{n},Ax_{n})+1)}{1+d(Bx_{n},Bx_{n+2})}\}
\\&&\nonumber = \max\{d(z_{n-1},z_{n+1}),\frac{d(z_{n-1},z_{n})(1+d(z_{n+1},z_{n+2}))}{1+d(z_{n-1},z_{n+1})},
\frac{ d(z_{n+1},z_{n+2})(d(z_{n-1},z_{n})+1)}{1+d(z_{n-1},z_{n+1})}\}
\\&&\nonumber = d(z_{n-1},z_{n+1}).
\end{eqnarray}
Hence
\begin{eqnarray}
 & d(z_{n},z_{n+2})&\leq \phi(d(z_{n-1},z_{n+1}))< d(z_{n-1},z_{n+1})
\end{eqnarray}
Thus the sequence $\{d(z_{n},z_{n+2})\}$ is decreasing and hence, is bounded. If, for some $n\in N,$ we have $d(z_{n-1},z_{n+1})\leq L $ and $d(z_{n},z_{n+2})> L,$ then from (8) we get
\begin{eqnarray}
\nonumber & d(z_{n},z_{n+2})&< L,
\end{eqnarray}
which is a contradiction. Then $\{d(z_{n},z_{n+2})\}$ is bounded. Now, if
\begin{eqnarray}
\lim\limits_{n\to\infty}d(z_{n},z_{n+2})=0
\end{eqnarray}
dose not hold, then there exists a subsequence $\{z_{n_{k}}\}$ of $\{z_{n}\}$ such that $\lim\limits_{n\to\infty}d(z_{n_k},z_{n_{k}+2})=s.$ From
$$
d(z_{n_{k}-1},z_{n_{k}+1})\leq d(z_{n_{k}-1},z_{n_{k}})+ d(z_{n_{k}},z_{n_{k}+2})+d(z_{n_{k}+1},z_{n_{k}+2})
$$
and
$$
d(z_{n_{k}},z_{n_{k}+2})\leq d(z_{n_{k}-1},z_{n_{k}})+ d(z_{n_{k}-1},z_{n_{k}+1})+d(z_{n_{k}+1},z_{n_{k}+2})
$$
we obtain that
$$
\lim\limits_{k\to +\infty}d(z_{n_{k}-1},z_{n_{k}+1})= s.
$$
Now, by (3) with $x = x_{n_{k}}$ and $y = x_{n_{k}+2},$ we have
\begin{eqnarray}
 & d(z_{n_{k}},z_{n_{k}+2})&\leq \phi(d(z_{n_{k}-1},z_{n_{k}+1}))
\end{eqnarray}
From (9) as $k\to\infty,$ we get $s\leq \phi(s)$ which implies $s = 0.$

Now, if possible, let $\{y_{n}\}$ be not a Cauchy sequence. Then there exists $\epsilon > 0 $ for which we can find subsequences $\{z_{n_{k}}\}$ and $\{z_{m_{k}}\}$ of $\{z_{n}\}$ with $n_{k}> m_{k}\geq k $ such that
\begin{eqnarray}
 & d(z_{n_{k}},z_{m_{k}})&\geq \epsilon.
\end{eqnarray}
Further, corresponding to $m_{k},$ we can choose $n_{k}$ in such a way that it is the smallest integer with $n_{k}-m_{k}\geq 4 $ and satisfying (11). Then
\begin{eqnarray}
 & d(z_{n_{k}-2},z_{m_{k}})& <  \epsilon.
\end{eqnarray}
Now, using (11), (12) and the rectangular inequality, we get
\begin{eqnarray}
 & \nonumber \epsilon \leq d(z_{n_{k}},z_{m_{k}})& \leq d(z_{n_{k}},z_{n_{k}-2})+d(z_{n_{k}-2},z_{n_{k}-1})+d(z_{n_{k}-1},z_{m_{k}})
 \\&& \nonumber < d(z_{n_{k}},z_{n_{k}-2})+d(z_{n_{k}-2},z_{n_{k}-1})+\epsilon.
\end{eqnarray}
Letting $k\to +\infty $ in the above inequality, using (6) and (9), we obtain
\begin{eqnarray}
\lim\limits_{k\to \infty }d(z_{n_{k}},z_{m_{k}})=\epsilon^{+}.
\end{eqnarray}
From
\begin{eqnarray}
\nonumber & &d(z_{n_{k}},z_{m_{k}})-d(z_{m_{k}},z_{m_{k}-1}) - d(z_{n_{k}-1},z_{n_{k}})
\\&& \nonumber \leq d(z_{n_{k}-1},z_{m_{k}-1})\leq d(z_{n_{k}-1},z_{n_{k}})+d(z_{m_{k}},z_{n_{k}})+d(z_{m_{k}-1},z_{m_{k}}),
\end{eqnarray}
letting $k\to +\infty,$ we obtain
\begin{eqnarray}
\lim\limits_{k\to\infty}d(z_{n_{k}-1},z_{m_{k}-1})=\epsilon.
\end{eqnarray}
From (3) with $ x = x_{n_{k}}$ and $y = x_{m_{k}},$  we get
\begin{eqnarray}
\nonumber &d(Ax_{m_{k}},Ax_{n_{k}})&\leq \phi(M(x_{m_{k}},x_{n_{k}}))+ C \min\{d(Bx_{n_{k}},Ax_{n_{k}}),d(Bx_{m_{k}},Ax_{m_{k}}),d(Bx_{n_{k}},Ax_{m_{k}}),d(Bx_{m_{k}},Ax_{n_{k}})\}
 \\&& \nonumber = \phi(M(x_{m_{k}},x_{n_{k}}))+ C \min\{d(z_{n_{k}-1},z_{n_{k}}),d(z_{m_{k}-1},z_{m_{k}}),d(z_{n_{k}-1},z_{m_{k}}),d(z_{m_{k}-1},z_{n_{k}})\}
\end{eqnarray}
where
\begin{eqnarray}
\nonumber &M(x_{m_{k}},x_{n_{k}})&=\max\{d(Bx_{m_{k}},Bx_{n_{k}}),\frac{ d(Bx_{m_{k}},Ax_{m_{k}})(d(Bx_{n_{k}},Ax_{n_{k}})+1)}{1+d(Bx_{m_{k}},Bx_{n_{k}})},
\frac{ d(Bx_{n_{k}},Ax_{n_{k}})(d(Bx_{m_{k}},Ax_{m_{k}})+1)}{1+d(Bx_{m_{k}},Bx_{n_{k}})}\}
\\&&\nonumber = \max\{d(z_{m_{k}-1},z_{n_{k}-1}),\frac{ d(z_{m_{k}-1},z_{m_{k}})(d(z_{n_{k}-1},z_{n_{k}})+1)}{1+d(z_{m_{k}-1},z_{n_{k}-1})},
\frac{ d(z_{n_{k}-1},z_{n_{k}})(d(z_{m_{k}-1},z_{m_{k}})+1)}{1+d(z_{m_{k}-1},z_{n_{k}-1})}\}.
\end{eqnarray}
Now, using the continuity of $\phi$ as $k\to +\infty,$ we obtain
 \begin{eqnarray}
\nonumber &\epsilon &\leq \phi(\epsilon)+0 < \epsilon,
\end{eqnarray}
which implies that $\epsilon = 0,$ a contradiction with $\epsilon > 0.$ Hence, $\{z_{n}\}$ is a G.M.S Cauchy sequence. Since $(BX,d)$ is complete G.M.S, there exists $z\in BX $ such that $\lim\limits_{n\to\infty}z_{n}= z.$ Let $u\in X$ be such that $Bu = z,$ applying (3) with $x = x_{n_{k}}$
\begin{eqnarray}
\nonumber &d(Au,Ax_{n_{k}})&\leq \phi(M(u,x_{n_{k}}))+L \min\{d(Bx_{n_{k}},Ax_{n_{k}}),d(Bu,Au),d(Bx_{n_{k}},Au),d(Bu,Ax_{n_{k}})\}
\\&& \nonumber = \phi(M(u,x_{n_{k}}))+L \min\{d(z_{n_{k}-1},z_{n_{k}}),d(Bu,Au),d(z_{n_{k}-1},Au),d(z,z_{n_{k}})\}
\\&& = \phi(M(u,x_{n_{k}}))+0,
\end{eqnarray}
where
\begin{eqnarray}
\nonumber &M(u,x_{n_{k}})&=\max\{d(Bu,Bx_{n_{k}}),\frac{d(Bu,Au)(d(Bx_{n_{k}},Ax_{n_{k}})+1)}{1+d(Bu,Bx_{n_{k}})},
\frac{ d(Bx_{n_{k}},Ax_{n_{k}})(d(Bu,Au)+1)}{1+d(Bu,Bx_{n_{k}})}\}
\\&&\nonumber = \max\{d(z,z_{n_{k}-1}),\frac{ d(Bu,Au)(d(z_{n_{k}-1},z_{n_{k}})+1)}{1+d(Bu,z_{n_{k}-1})},
\frac{ d(z_{n_{k}-1},z_{n_{k}})(d(Bu,Au)+1)}{1+d(Bu,z_{n_{k}-1})}\}
\\&&\nonumber = d(Bu,Au) \; as \; k\to \infty.
\end{eqnarray}
We get from (15) that
\begin{eqnarray}
\nonumber & d(Bu,Au)&\leq \liminf_{k\to\infty}[d(Bu,z_{n_{k}-1})+ d(z_{n_{k}-1},z_{n_{k}})+d(Au,Ax_{n_{k}})]\leq \liminf_{k\to\infty} d(Au,Ax_{n_{k}})
\\&& = \phi(d(Bu,Au))< d(Bu,Au),
 \end{eqnarray}
which implies $d(Bu,Au)=0,$ that is, $z = Bu = Au$ and so $z$ is a point coincidence for $A$ and $B.$

Now, we prove that $z$ is the unique point of coincidence  of $A$ and $B.$ Let $x$ and $y$ be arbitrary points coincidence of $A$ and $B$ such that $x = Au = Bu$ and $y = Av = Bv.$  Using the condition (3), it follows that
   \begin{eqnarray}
\nonumber &d(x,y)=d(Au,Av)&\leq \phi(\max\{d(Bu,Bv),\frac{d(Bu,Au)(d(Bv,Av)+1)}{1+d(Bu,Bv)},
\frac{ d(Bv,Av)(d(Bu,Au)+1)}{1+d(Bu,Bv)}\})
\\&&\nonumber + C \min\{d(Bv,Av),d(Bu,Au),d(Bv,Au),d(Bu,Av)\}
\\&& \nonumber = \phi(d(Bu,Bv))< d(Bu,Bv)=d(x,y),
\end{eqnarray}
which implies that $d(x,y)=0.$ Thus, $x = y$ and $A, B$ have a unique point of coincidence.

Next, we prove that $z = Az = Bz.$ If $z$ is a point of coincidence of $A$ and $B$ as $A$ and $B$ weakly compatible, we obtain that $Az = AAu =ABu = BAu =Bz$ and so $z = Az = Bz.$ Consequently, $z$ is unique common fixed point of $A$ and $B.$
\\
{\bf Example 3.1.} Suppose that $(X,d)$ as in example 2.1. let $A,B: X\to X $ and $\phi(t): [0,\infty)\to [0,\infty)$ defined by
\\
 $Ax =\frac{1}{2}x,\; Bx = x$
 and  $\phi(t)=\frac{t}{2},\; \forall \; t\in [0,\infty).$\\
 Then $A,B$ and $\phi$ satisfy all the conditions of Theorem 3.1. Hence, $0$  unique coincidence and common fixed point of $A$ and $B.$
\\
{\bf Corollary 3.1.} Replacing  the condition (3) in Theorem 3.1 with the following condition:
 \begin{eqnarray}
&\nonumber d(Ax,Ay)&\leq a_{1}d(Bx,By)+a_{2}\frac{ d(Bx,Ax)(d(By,Ay)+1)}{1+d(Bx,By)}+a_{3}\frac{ d(By,Ay)(d(Bx,Ax)+1)}{1+d(Bx,By)}
\\&&+ L\min\{d(Bx,Ax), d(By,Ay),d(Bx,Ay),d(By,Ax)\},
\end{eqnarray}
where $a_{1},a_{2},a_{3}, L\geq 0,$ and $a_{1}+a_{2}+a_{3}< 1.$
\\
{\bf Corollary 3.2.}  Putting $B =I$(the identity mapping) in Theorem 3.1. Then one can get a unique fixed point of $A.$
\\
{\bf Remark 3.1.} \cite[Theorem 7]{ARS2015} is spatial case of Theorem 3.1.

Next, we introduce some coincidence point theorems for two  $(\alpha,\psi,\phi)$-contractions self-mappings of rational type in complete G.M.S.
\\
{\bf Theorem 3.2.} Let $(X,d)$ be a G.M.S and  let $A,B: X\to X $ be two self-mappings satisfy the following conditions:
 \begin{eqnarray}
\phi(\beta(Bx,By)d(Ax,Ay))\leq \phi(M(x,y))-\psi(M(x,y)) \; \forall \; x,y \in X,
\end{eqnarray}
where $M(x,y)$ as in Theorem 3.1,  $AX\subset BX,$ and  $(BX,d)$ is a complete G.M.S
\\
Consider also that the next conditions hold:

(i) $\exists$ $x_{0}\in X $ such that $\beta(fx_{0},Tx_{0})\geq 1,$

(ii) $A$ is $B-\beta-$admissible,

(iii) $X$ is $\beta-$regular and   $\beta(x_{m},x_{n})\geq 1,$ for each $x_{n}\in X,$ and  $\forall\; m,n \in N, m\ne n,$

(iv) either $\beta(Bx,By)\geq 1$ or $\beta(By,Bx)\geq 1$ whenever $Bx = Ax$ and $By = Ay,$

(iiv) $\psi: [0,\infty) \to  [0,\infty) $ be a lower semi-continuous function and $\psi(t)=0 \; \iff t =0.$
\\
Then $A$ and $B$ have a unique point of coincidence in $X.$ Moreover if $A$ and $B$ are weakly compatible, then $A$ and $B$ have a unique common fixed point.
\\
{\bf Proof.} Suppose that $x_{0}\in X,$ $\beta(Bx_{0},Ax_{0})\geq 1.$  Define $\{z_n\}$ and $\{x_n\}$ be two sequences in $X$  such that $
z_n =Bx_{n+1}=Ax_{n}, \; n=0,1,2,3,....$
If $z_{n}=z_{n+1},$ then $Bx_{n+1}=Ax_{n+1}$  which implies that $x_{n+1}$ is a coincidence point of $A$ and $B.$ Consequently, we can suppose that $z_{n}\ne z_{n+1} $ for all $n\in N.$ From (i), we get that $\beta(Bx_{0},Ax_{0})=\beta(Bx_{0},Bx_{1})\geq 1.$ Also, by (ii) we have that $\beta(Ax_{0},Ax_{1})=\beta(Bx_{1},Bx_{2})\geq 1,$ $\beta(Ax_{1},Ax_{2})=\beta(Bx_{2},Bx_{3})\geq 1.$ continuous with this process we obtain that  $\beta(Bx_{n},Bx_{n+1})\geq 1.$ Now, by using (18), we get
\begin{eqnarray}
 & \phi(d(Ax_{n},Ax_{n+1}))&\leq \phi(\beta(Bx_{n},Bx_{n+1})d(Ax_{n},Ax_{n+1}))\leq \phi(M(x_{n},x_{n+1}))-\psi(M(x_{n},x_{n+1}))
\end{eqnarray}
where
\begin{eqnarray}
\nonumber &M(x_{n},x_{n+1})&=\max\{d(Bx_{n},Bx_{n+1}),\frac{ d(Bx_{n},Ax_{n})(d(Bx_{n+1},Ax_{n+1})+1)}{1+d(Bx_{n},Bx_{n+1})}, \frac{ d(Bx_{n+1},Ax_{n+1})(d(Bx_{n},Ax_{n})+1)}{1+d(Bx_{n},Bx_{n+1})}\}
\\&&\nonumber = \max\{d(z_{n-1},z_{n}),\frac{d(z_{n-1},z_{n})(1+d(z_{n},z_{n+1}))}{1+d(z_{n-1},z_{n})}, d(z_{n},z_{n+1})\},
\end{eqnarray}
we consider the following cases

$\bullet$ If $M(x_{n},x_{n+1})= d(z_{n-1},z_{n})$ from $(20)$ we have
\begin{eqnarray}
 &\nonumber  \phi(d(z_{n},z_{n+1}))&\leq \phi(d(z_{n-1},z_{n}))- \psi(d(z_{n-1},z_{n})) < \phi(d(z_{n-1},z_{n})),
\end{eqnarray}
Since $\phi$ is nondecreasing we have
\begin{eqnarray}
 d(z_{n},z_{n+1})< \phi(d(z_{n-1},z_{n})).
\end{eqnarray}

 $\bullet$ If $M(x_{n},x_{n+1})= \frac{d(z_{n-1},z_{n})(1+d(z_{n},z_{n+1}))}{1+d(z_{n-1},z_{n})}$ from $(20)$ we obtain
 \begin{eqnarray}
 \nonumber & \phi(d(z_{n},z_{n+1}))&\leq \phi(\frac{d(z_{n-1},z_{n})(1+d(z_{n},z_{n+1}))}{1+d(z_{n-1},z_{n})})-\psi(\frac{d(z_{n-1},z_{n})(1+d(z_{n},z_{n+1}))}{1+d(z_{n-1},z_{n})})
  \\&&\nonumber  < \phi(\frac{d(z_{n-1},z_{n})(1+d(z_{n},z_{n+1}))}{1+d(z_{n-1},z_{n})}).
\end{eqnarray}
The nondecreasing property of $ \phi $ implies that
 \begin{eqnarray}
 \nonumber & &d(z_{n},z_{n+1}) < \frac{d(z_{n-1},z_{n})(1+d(z_{n},z_{n+1}))}{1+d(z_{n-1},z_{n})}
 \\&& \nonumber  \Longrightarrow d(z_{n},z_{n+1})+d(z_{n},z_{n+1})d(z_{n-1},z_{n})< d(z_{n-1},z_{n})+d(z_{n},z_{n+1})d(z_{n-1},z_{n})
 \\&&   \Longrightarrow d(z_{n},z_{n+1}) < d(z_{n-1},z_{n}).
\end{eqnarray}
Hence, (20) is obtained.

$\bullet$ If $M(x_{n},x_{n+1})= d(z_{n},z_{n+1})).$ By $(18)$ we obtain
\begin{eqnarray}
\nonumber  & \phi(d(z_{n},z_{n+1}))& \leq \phi(d(z_{n},z_{n+1}))-\psi(d(z_{n},z_{n+1}))
\\&&\nonumber < \phi(d(z_{n},z_{n+1})),
\end{eqnarray}
this is a contradiction.

In any case, we proved that (18) holds. since $\{d(z_{n},z_{n+1})\}$ is decreasing. Hence, it converges to a nonnegative number, $s\geq 0.$ If $s > 0,$ then letting $n\to +\infty$ in (18), we deduce
\begin{eqnarray}
\nonumber & s &\leq \phi(\max\{s,\frac{s(1+s)}{1+s}, s\})= \phi(s)< s,
\end{eqnarray}
which implies that $s = 0,$ that is
\begin{eqnarray}
& \lim\limits_{n\to\infty}d(z_{n},z_{n+1})& = 0.
\end{eqnarray}
Suppose that $z_{n}\ne z_{m}$ for all $m\ne n$ and prove that $\{z_{n}\}$ is G.M.S Cauchy sequence. First, we show that the sequence $\{d(z_{n},z_{n+2})\}$ is bounded. Since $\lim\limits_{n\to\infty}d(z_{n},z_{n+1})=0,$ there exists $L > 0 $ such that $d(z_{n},z_{n+1})\leq L $ for all $n\in N.$ If $d(z_{n},z_{n+2})> L $ for all $n\in N,$ from (18) we have
\begin{eqnarray}
\nonumber & \phi(d(z_{n},z_{n+2}))&\leq \phi(\beta(d(Bx_{n},Bx_{n+2})d(Ax_{n},Ax_{n+2}))=d(Ax_{n},Ax_{n+2})
\\&&\nonumber \leq \phi(M(x_{n},x_{n+2}))-\psi(M(x_{n},x_{n+2}))
\\&&<\phi(M(x_{n},x_{n+1}))= \phi(d(z_{n-1},z_{n+1}))\;\; as \;\; n\to \infty,
\end{eqnarray}
where
\begin{eqnarray}
\nonumber &M(x_{n},x_{n+2})&=\max\{d(Bx_{n},Bx_{n+2}),\frac{ d(Bx_{n},Ax_{n})(d(Bx_{n+2},Ax_{n+2})+1)}{1+d(Bx_{n},Bx_{n+2})},
\frac{ d(Bx_{n+2},Ax_{n+2})(d(Bx_{n},Ax_{n})+1)}{1+d(Bx_{n},Bx_{n+2})}\}
\\&&\nonumber = \max\{d(z_{n-1},z_{n+1}),\frac{d(z_{n-1},z_{n})(1+d(z_{n+1},z_{n+2}))}{1+d(z_{n-1},z_{n+1})},
\frac{ d(z_{n+1},z_{n+2})(d(z_{n-1},z_{n})+1)}{1+d(z_{n-1},z_{n+1})}\}
\\&&\nonumber = d(z_{n-1},z_{n+1})  \;\;\textrm{as}\;\; n \to \infty.
\end{eqnarray}
Hence
\begin{eqnarray}
 & d(z_{n},z_{n+2})&\leq \phi(d(z_{n-1},z_{n+1}))< d(z_{n-1},z_{n+1})
\end{eqnarray}
Thus the sequence $\{d(z_{n},z_{n+2})\}$ is decreasing and hence, is bounded. If, for some $n\in N,$ we have $d(z_{n-1},z_{n+1})\leq L $ and $d(z_{n},z_{n+2})> L,$ then from (24) we get
\begin{eqnarray}
\nonumber & d(z_{n},z_{n+2})&< L,
\end{eqnarray}
which is a contradiction. Then $\{d(z_{n},z_{n+2})\}$ is bounded. Now, if
\begin{eqnarray}
\lim\limits_{n\to\infty}d(z_{n},z_{n+2})=0
\end{eqnarray}
dose not hold, then there exists a subsequence $\{z_{n_{k}}\}$ of $\{z_{n}\}$ such that $\lim\limits_{n\to\infty}d(z_{n_k},z_{n_{k}+2})=s.$ From
$$
d(z_{n_{k}-1},z_{n_{k}+1})\leq d(z_{n_{k}-1},z_{n_{k}})+ d(z_{n_{k}},z_{n_{k}+2})+d(z_{n_{k}+1},z_{n_{k}+2})
$$
and
$$
d(z_{n_{k}},z_{n_{k}+2})\leq d(z_{n_{k}-1},z_{n_{k}})+ d(z_{n_{k}-1},z_{n_{k}+1})+d(z_{n_{k}+1},z_{n_{k}+2})
$$
we obtain that
$$
\lim\limits_{k\to +\infty}d(z_{n_{k}-1},z_{n_{k}+1})= s.
$$
Now, by (18) one can obtain that
\begin{eqnarray}
 &\phi( d(z_{n_{k}},z_{n_{k}+2}))& \leq  \phi(d(z_{n_{k}-1},z_{n_{k}+1}))- \psi(d(z_{n_{k}-1},z_{n_{k}+1}))
 \Longrightarrow \phi( s)  <  \phi(s) \; \textrm{as} \; n\to \infty.
\end{eqnarray}
 which implies $s = 0.$

Now, if possible, let $\{z_{n}\}$ be not a Cauchy sequence. Then there exists $\epsilon > 0 $ for which we can find subsequences $\{z_{n_{k}}\}$ and $\{z_{m_{k}}\}$ of $\{z_{n}\}$ with $n_{k}> m_{k}\geq k $ such that
\begin{eqnarray}
 & d(z_{n_{k}},z_{m_{k}})&\geq \epsilon.
\end{eqnarray}
Further, corresponding to $m_{k},$ we can choose $n_{k}$ in such a way that it is the smallest integer with $n_{k}-m_{k}\geq 4 $ and satisfying (27). Then
\begin{eqnarray}
 & d(z_{n_{k}-2},z_{m_{k}})& <  \epsilon.
\end{eqnarray}
Now, using (27), (28) and the rectangular inequality, we get
\begin{eqnarray}
 & \nonumber \epsilon \leq d(z_{n_{k}},z_{m_{k}})& \leq d(z_{n_{k}},z_{n_{k}-2})+d(z_{n_{k}-2},z_{n_{k}-1})+d(z_{n_{k}-1},z_{m_{k}})
 \\&& \nonumber < d(z_{n_{k}},z_{n_{k}-2})+d(z_{n_{k}-2},z_{n_{k}-1})+\epsilon.
\end{eqnarray}
Letting $k\to +\infty $ in the above inequality, using (22) and (25), we obtain
\begin{eqnarray}
\lim\limits_{k\to \infty }d(z_{n_{k}},z_{m_{k}})=\epsilon^{+}.
\end{eqnarray}
From
\begin{eqnarray}
\nonumber & &d(z_{n_{k}},z_{m_{k}})-d(z_{m_{k}},z_{m_{k}-1}) - d(z_{n_{k}-1},z_{n_{k}})
\\&& \nonumber \leq d(z_{n_{k}-1},z_{m_{k}-1})\leq d(z_{n_{k}-1},z_{n_{k}})+d(z_{m_{k}},z_{n_{k}})+d(z_{m_{k}-1},z_{m_{k}}),
\end{eqnarray}
letting $k\to +\infty,$ we obtain
\begin{eqnarray}
\lim\limits_{k\to\infty}d(z_{n_{k}-1},z_{m_{k}-1})=\epsilon.
\end{eqnarray}
From (18) with $ x = x_{n_{k}}$ and $y = x_{m_{k}},$  we get
\begin{eqnarray}
\nonumber &\phi(d(Ax_{m_{k}},Ax_{n_{k}}))\leq \phi(\beta(d(Bx_{m_{k}},Bx_{n_{k}}))d(Ax_{m_{k}},Ax_{n_{k}}))&\leq \phi(M(x_{m_{k}},x_{n_{k}}))- \psi(M(x_{m_{k}},x_{n_{k}}))
\end{eqnarray}
where
\begin{eqnarray}
\nonumber &M(x_{m_{k}},x_{n_{k}})&=\max\{d(Bx_{m_{k}},Bx_{n_{k}}),\frac{ d(Bx_{m_{k}},Ax_{m_{k}})(d(Bx_{n_{k}},Ax_{n_{k}})+1)}{1+d(Bx_{m_{k}},Bx_{n_{k}})},
\frac{ d(Bx_{n_{k}},Ax_{n_{k}})(d(Bx_{m_{k}},Ax_{m_{k}})+1)}{1+d(Bx_{m_{k}},Bx_{n_{k}})}\}
\\&&\nonumber = \max\{d(z_{m_{k}-1},z_{n_{k}-1}),\frac{ d(z_{m_{k}-1},z_{m_{k}})(d(z_{n_{k}-1},z_{n_{k}})+1)}{1+d(z_{m_{k}-1},z_{n_{k}-1})},
\frac{ d(z_{n_{k}-1},z_{n_{k}})(d(z_{m_{k}-1},z_{m_{k}})+1)}{1+d(z_{m_{k}-1},z_{n_{k}-1})}\}.
\end{eqnarray}
Now, using the continuity of $\phi$ as $k\to +\infty,$ we obtain
 \begin{eqnarray}
\nonumber &\epsilon &\leq \phi(\epsilon) < \epsilon,
\end{eqnarray}
a contradiction is obtained with $\epsilon > 0,$ then $\epsilon = 0,$   Hence, $\{z_{n}\}$ is a G.M.S Cauchy sequence. Since $(BX,d)$ is  complete G.M.S , there exists $z\in BX $ such that $\lim\limits_{n\to\infty}z_{n}= z.$ Let $w\in X$ be such that $Bu = z,$ applying (18) with $x = x_{n_{k}}$
\begin{eqnarray}
\nonumber &\phi(d(Au,Ax_{n_{k}}))&\leq \phi(M(u,x_{n_{k}}))- \phi(M(u,x_{n_{k}})),
\end{eqnarray}
where
\begin{eqnarray}
\nonumber &M(u,x_{n_{k}})&=\max\{d(Bu,Bx_{n_{k}}),\frac{d(Bu,Au)(d(Bx_{n_{k}},Ax_{n_{k}})+1)}{1+d(Bu,Bx_{n_{k}})},
\frac{ d(Bx_{n_{k}},Ax_{n_{k}})(d(Bu,Au)+1)}{1+d(Bu,Bx_{n_{k}})}\}
\\&&\nonumber = \max\{d(z,z_{n_{k}-1}),\frac{ d(Bu,Au)(d(z_{n_{k}-1},z_{n_{k}})+1)}{1+d(Bu,z_{n_{k}-1})},
\frac{ d(z_{n_{k}-1},z_{n_{k}})(d(Bu,Au)+1)}{1+d(Bu,z_{n_{k}-1})}\}
\\&&\nonumber = d(Bu,Au) \; as \; k\to \infty.
\end{eqnarray}
We get from (15) that
\begin{eqnarray}
\nonumber & d(Bu,Au)&\leq \liminf_{k\to\infty}[d(Bu,z_{n_{k}-1})+ d(z_{n_{k}-1},z_{n_{k}})+d(Au,Ax_{n_{k}})]\leq \liminf_{k\to\infty} d(Au,Ax_{n_{k}})
\\&& = \phi(d(Bu,Au))< d(Bu,Au),
 \end{eqnarray}
which implies $d(Bu,Au)=0,$ that is, $z = Bu = Au$ and so $z$ is a coincidence point for $A$ and $B.$

Now, we prove that $z$ is the unique coincidence point of $A$ and $B.$ Let $x$ and $y$ be arbitrary coincidence points of $A$ and $B$ such that $x = Au = Bu$ and $y = Av = Bv.$  Using the condition (18), it follows that
   \begin{eqnarray}
\nonumber &\phi(d(x,y))=\phi(d(Au,Av))&\leq \phi(\max\{d(Bu,Bv),\frac{d(Bu,Au)(d(Bv,Av)+1)}{1+d(Bu,Bv)},
\frac{ d(Bv,Av)(d(Bu,Au)+1)}{1+d(Bu,Bv)}\})
\\&&\nonumber -\psi(\max\{d(Bu,Bv),\frac{d(Bu,Au)(d(Bv,Av)+1)}{1+d(Bu,Bv)},
\frac{ d(Bv,Av)(d(Bu,Au)+1)}{1+d(Bu,Bv)}\})
\\&& \nonumber = \phi(d(Bu,Bv))- \psi(d(Bu,Bv)) < \phi(d(Bu,Bv))=\phi(d(x,y)),
\end{eqnarray}
which implies that $d(x,y)=0.$ Thus, $x = y$ and $A, B$ have a unique coincidence point.

As in the conclusion in last paragraph of the proof of Theorem 3.1  and the weakly compatible property of $A$ and $B,$ we obtain that $A$ and $B$ is unique common fixed point.
\\
{\bf Corollary 3.3.}  Putting $B =I$ in Theorem 3.2. Then one can get a unique fixed point of $A.$
\\
{\bf Remark 3.2.} \cite[Theorem 16]{ARS2015} is spatial case of Theorem 3.2.
\section{An application in dynamical programming }
The aim of this section is to use Theorem 3.1 to study the existence and uniqueness of solutions of the following system of functional equations:
\begin{equation}
\begin{array}{lll}
  w(x)=\sup\limits_{b\in E}\{h(a,b)+F(a,b,z(G(a,b))), \\
z(x)=\sup\limits_{b\in E}\{h(a,b)+F(a,b,w(G(a,b)))\}, \\
\end{array}
\end{equation}
which use in dynamic programming (see \cite{liu2001,chang91,chang1991}), where $E$ is a state space, $S$ is a decision space and $x\in S,$ $y\in E,$  $h:S\times E\to \Re ,$ $G: S\times E \to S$ and $F: S\times E\times \Re \to \Re $ are considered operators.

We denote by $B(S)$ the set of all bounded functionals on $S.$ Also, we define $\parallel . \parallel_{\infty} $ by
$$
\parallel v\parallel_{\infty}=\sup\limits_{x\in S}\mid v(x) \mid,  \; \forall \; v\in B(S).
$$
{\bf Remark}. we note that the space $(B(S), \parallel . \parallel_{\infty} )$ is a Banach, where the distance function in $B(S)$ defined as follows:
$$
d_{\infty}(T_{1},T_{2})=\sup\limits_{x\in S} \mid T_{1}(x)- T_{2}(x) \mid \; \forall \; T_{1},T_{2}\in B(S).
$$
\\
{\bf Lemma 4.1} (\cite{ARS2015}). Let $F_{1}, F_{2} : S \to \Re $ are bounded functionals, then
\begin{eqnarray}
\mid \sup\limits_{x\in S}F_{1}(x)- \sup\limits_{x\in S}F_{2}(x)\mid\leq \sup\limits_{x\in S}\mid F_{1}(x)- F_{1}(x)\mid.
\end{eqnarray}
{\bf Proposition 4.1}. Suppose that $h, F(.,.,0),F(.,.,1): S\times E \to \Re $ are three bounded functionals and
\begin{equation}
\mid F(a,b,t_{1}) - F(a,b,t_{2})\mid \leq C \mid t_{1} - t_{2} \mid, C\geq 0, \; \forall \; a\in S,\; b\in E \;\textrm{and} \; t_{1},t_{2}\in \Re.
\end{equation}
Also, let $O : B(S)\to B(S)$ be an operator defined as follows:
\begin{equation}
\begin{array}{lll}
  (O w)(b)=\sup\limits_{b\in E}\{h(a,b)+F(a,b,z(G(a,b))\}, \; \forall \; a\in S,\\
z(a)=\sup\limits_{b\in E}\{h(a,b)+F(a,b,w(G(a,b)))\}, \; \forall \; a\in S.\\
\end{array}
\end{equation}
For all $w\in B(S)$ and $x\in S.$
Then $O$ is well defined.
\\
{\bf Proof.} It is enough to prove that $O : B(S)\to B(S)$ is bounded for all $w\in B(S).$ From the boundedness of  $w,$ $h$ and $F$ we have
\begin{eqnarray}
\nonumber &|  (O w)(x) | &\leq \sup\limits_{b\in E}|h(a,b)+F(a,b,z(G(a,b))|
\\&&\nonumber \leq  \sup\limits_{b\in E}|h(a,b)|+\sup\limits_{b\in E}|F(a,b,z(G(a,b))-F(a,b,0)| + \sup\limits_{b\in E} |F(a,b,0)-F(a,b,1)|+\sup\limits_{b\in E} |F(a,b,1)|
\\&&\nonumber  \leq \sup\limits_{b\in E}|h(a,b)|+\sup\limits_{b\in E}|z(G(a,b))| +1+ \sup\limits_{b\in E} |F(a,b,1)|
\\&& \nonumber \leq \sup\limits_{b\in E}|h(a,b)|+\sup\limits_{b\in E}|\sup\limits_{b\in E}\{h(a,b)+F(a,b,w(G(a,b)))\}| + 1+ \sup\limits_{b\in E} |F(a,b,1)|
\\&& \leq 2\sup\limits_{b\in E}|h(a,b)|+\sup\limits_{b\in E}|w(G(a,b))| + 2 + 2\sup\limits_{b\in E} |F(a,b,1)|\leq C^{'},\;\; C^{'}> 0.
\end{eqnarray}
 Which give the boundedness of $O(w)$ on $w.$ Hence, $R$ is well defined.
 \\
 {\bf Theorem 4.1}. Consider the assumptions of  Proposition 4.1 and the following property:
 \begin{eqnarray}
&d(F(a,b,z(w_{1}(G(a,b))),F(a,b,z(w_{2}(G(a,b))))&\leq \phi(M(w_{1},w_{2})) + C m(w_{1},w_{2})
\end{eqnarray}
   where
 \begin{eqnarray}
&&\nonumber M(w_{1},w_{2})=\max\{d_{\infty}(zw_{1},zw_{2}),\frac{ d_{\infty}(zw_{1},Ow_{1})(d_{\infty}(zw_{2},Ow_{2})+1)}{1+d_{\infty}(zw_{1},zw_{2})}, \frac{ d_{\infty}(zw_{2},Ow_{2})(d_{\infty}(zw_{1},Ox)+1)}{1+d_{\infty}(zw_{1},zw_{2})}\}
\\&&\nonumber  m(w_{1},w_{2})= \min\{d_{\infty}(zw_{1},Ow_{1}), d_{\infty}(zw_{2},Ow_{2}),d_{\infty}(zw_{1},Ow_{2}),d_{\infty}(zw_{2},Ow_{1})\}.
\end{eqnarray}
For all $w_{1}, w_{2}\in B(S),$ and all $a\in S,$ all $b\in E.$ Also, the function $\phi$ as in Theorem 3.1.
\\
 Then (32) has a unique common solution $w_{0}\in B(S).$
 \\
 {\bf Proof. } First, we prove that the mappings in system (35) satisfy the condition (3). Indeed, by using Lemma 4.1, we have that, $\forall$ $w_{1}, w_{2}\in B(S),$ $\forall$ $x\in S,$
 \begin{eqnarray}
&&\nonumber d_{\infty}(Ow_{1},Ow_{2})\leq \sup\limits_{b\in E}\mid F(a,b,z(w_{1})-F(a,b,z(w_{2})\mid
\\&& \nonumber \leq \phi(M(w_{1},w_{2})) + C m(w_{1},w_{2}).
\end{eqnarray}
then all the conditions of Theorem 3.1 are satisfied, hence  the system (32) has a unique solution.
\\
{\bf Competing interests}
\\
The authors declare that they have no competing interests.
\\
{\bf Authors' contributions}
\\
The two authors contributed equally to this work. Both authors read and approved the final manuscript.
 
\end{document}